\documentclass[10pt]{article}
\usepackage{amsmath, amsthm, amssymb, amsfonts, mathrsfs, mathtools}
\usepackage{graphicx}
\usepackage{makeidx}
\usepackage[left=2cm,right=2cm,top=2cm,bottom=2cm]{geometry}
\usepackage{caption}
\usepackage{subcaption}
\usepackage[section]{placeins}
\usepackage{enumitem}
\usepackage{tikz}
\usepackage{stmaryrd}

\usetikzlibrary{decorations.pathreplacing}
\tikzstyle{vertex}=[auto=left,circle,draw=black,fill=white, inner sep=1.5]

\usepackage{hyperref}
\hypersetup{colorlinks=true, citecolor={blue}, linkcolor=blue, filecolor=magenta, urlcolor=cyan}

\newtheorem{theorem}{Theorem}[section]

\newtheorem{lema}[theorem]{Lemma}
\newtheorem{corollary}{Corollary}[theorem]

\newtheorem{rem}{Remark}[section]

\title{Pretty good fractional revival on abelian Cayley graphs}

\author{ Akash Kalita and Bikash Bhattacharjya\\
Department of Mathematics\\
Indian Institute of Technology Guwahati, India\\
akash.kalita@iitg.ac.in, b.bikash@iitg.ac.in }

\date{}
\begin{document}
\maketitle

\vspace{-0.3in}

\begin{center}{\textbf{Abstract}}\end{center}
Let $\Gamma$ be a graph with the adjacency matrix $A$. The transition matrix of $\Gamma$, denoted $H(t)$, is defined as $H(t) := \exp(-\textbf{i}tA)$, where $\textbf{i} := \sqrt{-1}$ and $t$ is a real variable. The graph $\Gamma$ is said to exhibit fractional revival (FR in short) between the vertices $a$ and $b$ if there exists a positive real number $t$ such that $H(t){\textbf{e}_{a}} = \alpha{\textbf{e}_{a}} + \beta{\textbf{e}_{b}}$, where $\alpha, \beta \in \mathbb{C}$ such that $\beta \neq 0$ and $|\alpha|^2 + |\beta|^2 = 1$. The graph $\Gamma$ is said to exhibit pretty good fractional revival (PGFR in short) between the vertices $a$ and $b$ if there exists a sequence of real numbers $\{t_k\}$ with $\lim_{k\to\infty} H(t_k){\textbf{e}_{a}} = \alpha{\textbf{e}_{a}} + \beta{\textbf{e}_{b}}$, where $\alpha, \beta \in \mathbb{C}$ such that $\beta \neq 0$ and $|\alpha|^2 + |\beta|^2 = 1$. In the definition of PGFR, if $\alpha=0$ then $\Gamma$ is said to exhibit pretty good state transfer (PGST in short) between $a$ and $b$. 
In this paper, we obtain some sufficient conditions for circulant graphs exhibiting PGFR. We also find some sufficient conditions for non-circulant abelian Cayley graphs exhibiting PGFR. From these 
sufficient conditions, we find infinite families of circulant graphs and non-circulant abelian Cayley graphs exhibiting PGFR that fail to exhibit FR and PGST. Finally, we obtain some necessary conditions for some families of circulant graphs exhibiting PGFR. Some of our results generalize the results of Chan et al. [Pretty good quantum fractional revival in paths and cycles. \textit {Algebr. Comb.} 4(6) (2021), 989-1004.] for cycles.

\vspace*{0.3cm}
\noindent 
\textbf{Keywords:} Pretty good state transfer, fractional revival, pretty good fractional revival, Cayley graph\\
\textbf{Mathematics Subject Classifications:} 15A16, 05C50, 81P45 

\section{Introduction}\label{sec 1}
Continuous-time quantum walk serves as a foundational method within the realm of quantum information. Farhi and Gutmann~\cite{Farhi} used continuous-time quantum walk to study quantum algorithmic problems on graphs in 1998. The generation of entanglement serves as an important resource in the field of quantum information. Fractional revival can be used for entanglement generation in quantum information. Fractional revival takes place when a continuous-time quantum walk maps the state of a vertex into a superposition of the states of a subset of vertices which includes the initial vertex. Since 1989, the concept of fractional revival was studied extensively in the context of physics, see for example,~\cite{Averbukh,Chen,Christandl,Genest,Robinett,FR First,Spanner}. Later in 2019, Chan et al.~\cite{Chan FR} studied fractional revival in the context of algebraic graph theory.

Throughout the paper, the word graph always means a simple undirected graph. Let $\Gamma$ be a graph on $n$ vertices with the vertex set $V(\Gamma)$. For $a \in V(\Gamma)$, the \textit{vertex state} of $a$ is defined to be the vector $\textbf{e}_{a}$, where $\{\textbf{e}_{a} \colon a \in V(\Gamma)\}$ is the standard basis of ${\mathbb{C}}^n$. If $A$ is the adjacency matrix of $\Gamma$, then the \textit{transition matrix} of $\Gamma$ is defined by
\begin{align*}
H(t) :=  \exp(-\mathbf{i}t A) = \sum_{j = 0}^{\infty} \frac{(-\mathbf{i}t A)^j}{j!},~\mathrm{where}~\mathbf{i} := \sqrt{-1}~\mathrm{and}~t~\mathrm{is~a~real~number}. 
\end{align*}
The graph $\Gamma$ exhibits \textit{fractional revival} (FR in short) between two distinct vertices $a$ and $b$ if there exists a positive real number $t$ such that 
$$H(t){\textbf{e}_{a}} = \alpha{\textbf{e}_{a}} + \beta{\textbf{e}_{b}},$$ 
where $\alpha, \beta \in \mathbb{C}$ such that $\beta \neq 0$ and $|\alpha|^2 + |\beta|^2 = 1$. In the definition of FR, if $|\alpha|=|\beta|$ then $\Gamma$ is said to exhibit \textit{balanced} FR between $a$ and $b$. In the definition of FR, if $\alpha = 0$ then $\Gamma$ is said to exhibit \textit{perfect state transfer} (PST in short) between $a$ and $b$. The graph $\Gamma$ is said to exhibit  \textit{periodicity} if there exists a positive real number $t$ and a complex number $\delta$ of unit modulus such that $H(t) = \delta I$.  

Bernard et al.~\cite{Bernard} provided some graph exhibiting balanced FR between antipodal vertices. Chan et al.~\cite{Chan FR} proved that a path on $n$ vertices exhibits FR if and only if $n \in \{2,3,4\}$. Chan et al.~\cite{Association Scheme} explored FR on graphs whose adjacency matrices belong to the Bose–Mesner algebra of association schemes. Godsil and Zhang~\cite{non-cospectral} provided classes of graphs exhibiting FR between non-cospectral vertices. Monterde~\cite{twin vertices} characterized the existence of FR between twin vertices on weighted graphs with respect to adjacency matrix, Laplacian matrix and signless Laplacian matrix. Wang et al.~\cite{J. Wang} and, Cao and Luo~\cite{Cao} studied FR on Cayley graphs over abelian groups. Cao and Luo~\cite{Cao} provided several new classes of Cayley graphs over abelian groups exhibiting FR. Wang et al.~\cite{Semi Cayley Graph} provided a characterization for the existence of FR on semi-Cayley graphs over abelian groups. Recently, Jitngam et al.~\cite{Finite Commutative Ring} classified all the finite local rings such that the corresponding unitary Cayley graphs exhibit FR. They also provided sufficient conditions for the existence of FR on unitary Cayley graphs over finite commutative rings.  

Chan et al.~\cite{Chan FR} developed tools that can be used to construct graphs exhibiting FR. They proved that a cycle exhibits FR if and only if the number of vertices of the cycle is either $4$ or $6$. Two vertices $a$ and $b$ of a graph $\Gamma$ are called \textit{cospectral}, if the subgraphs $\Gamma \setminus a$ and $\Gamma \setminus b$ hold the same list of eigenvalues. Chan et al.~\cite{Chan FR} proved that for a positive integer $k$, there are only finitely many connected graphs with maximum degree at most $k$ exhibiting FR between cospectral vertices. In the definition of FR, Chan et al.~\cite{Chan FR Laplacian} considered the Laplacian matrix in place of adjacency matrix to define Laplacian FR. They proved that Laplacian FR is a rare phenomenon. 

Chan et al.~\cite{PGFR Adjacency} introduced the notion of pretty good fractional revival, a relaxation to the notion of fractional revival on graphs. The graph $\Gamma$ is said to exhibit \textit{pretty good fractional revival} (PGFR in short) between two distinct vertices $a$ and $b$ if there exists a sequence of real numbers $\{t_k\}$ such that 
$$\lim_{k\to\infty} H(t_k){\textbf{e}_{a}} = \alpha{\textbf{e}_{a}} + \beta{\textbf{e}_{b}},$$ 
where $\alpha$ and $\beta$ are complex numbers with $\beta \neq 0$ and $|\alpha|^2 + |\beta|^2 = 1$. Observe that PGFR is an approximation to FR. Thus PGFR can be used for approximated entanglement generation in quantum information, whenever FR is absent in the corresponding network. For two distinct vertices $a$ and $b$ of $\Gamma$, we assume  that the first two rows and the first two columns of $A$ are indexed by $a$ and $b$, respectively. One can prove that $\Gamma$ exhibits PGFR between $a$ and $b$ if and only if there exists a sequence of real numbers $\{t_k\}$ such that 
\begin{align*}
\lim_{k\to\infty} H(t_k) =  \begin{pmatrix}
        M & \Large{\textbf{0}}  \\
       {\Large{\textbf{0}}}^T & N
       \end{pmatrix},
\end{align*}
where $M = \begin{pmatrix}
        \alpha & \beta  \\
       \beta & \gamma
       \end{pmatrix}, \alpha, \beta, \gamma \in \mathbb{C}$ with $\beta \neq 0$, $|\alpha|^2 + |\beta|^2 = 1~\mathrm{ and }~\gamma = -\frac{\overline \alpha \beta}{\overline \beta}$; $\Large{\textbf{0}}$ is the $2 \times (n - 2)$ zero matrix and $N$ is a $(n - 2) \times (n - 2)$ unitary matrix. In the definition of PGFR, if $\alpha = 0$ then $\Gamma$ is said to exhibit \textit{pretty good state transfer} (PGST in short) between $a$ and $b$. 

In~\cite{PGFR Adjacency}, Chan et al. developed the theory of PGFR on graphs with respect to the adjacency matrix and gave a complete classification of PGFR on paths and cycles. They proved that a cycle exhibits PGFR if and only if the number of vertices is of the form $2 p^s$, where $p$ is a prime number and $s$ is a positive integer, and it occurs between every pair of antipodal vertices. In~\cite{PGFR Laplacian}, Chan et al. studied Laplacian PGFR on graphs and classified paths and double stars that exhibit Laplacian PGFR. Drazen et al.~\cite{PGFR diagonal perturbation} generalized the concept of PGFR between two vertices of a graph to an arbitrary-sized subsets. Wang et al.~\cite{PGFR Nonabelian} provided necessary and sufficient conditions for the existence of PGFR on normal Cayley graphs over dicyclic groups. They used these conditions to obtain Cayley graphs exhibiting or not exhibiting PGFR. Pal and Bhattacharjya~\cite{Pal 2} proved that a cycle and its complement exhibit PGST if and only if the number of vertices is a power of two, and it occurs between antipodal vertices. Further, Pal~\cite{Pal, State transfer on circulant Pal} generalized this result to a class of circulant graphs. In this paper, we study the existence of PGFR on Cayley graphs over abelian groups.

The organization of the paper is as follows. In Section~\ref{sec 2}, we discuss some basic definitions and results which will be required in later sections. In Section~\ref{sec 3}, we determine a necessary and sufficient condition for the existence of PGFR on Cayley graphs over abelian groups. In Section~\ref{sec 4}, we obtain some sufficient conditions for circulant graphs exhibiting PGFR [see Theorem~\ref{PGFR existence example} and Theorem~\ref{PGFR first existence example}] using the necessary and sufficient condition described in Section~\ref{sec 3}. 
From these sufficient conditions, we find infinite families of circulant graphs exhibiting PGFR that fail to exhibit FR and PGST. We also obtain some necessary conditions for circulant graphs exhibiting PGFR [see Theorem~\ref{1st PGFR Nonexistence}, Theorem~\ref{Complement does not exhibit PGFR}, Theorem~\ref{2nd PGFR Nonexistence} and Theorem~\ref{another family not exhibiting PGFR}]. In Section~\ref{sec 5}, we obtain some sufficient conditions for non-circulant abelian Cayley graphs exhibiting PGFR. From these sufficient conditions, we find infinite families of non-circulant abelian Cayley graphs exhibiting PGFR that fail to exhibit FR and PGST. We also obtain an infinite family of non-circulant abelian Cayley graphs that fail to exhibit PGFR.

\section{Preliminaries}\label{sec 2}
Let $G$ be a finite group and $S$ be a subset of $G \setminus \{\textbf{1}\}$ with $S = \{y^{-1} \colon y \in S\}$, where $\textbf{1}$ is the identity element of $G$. The \textit{Cayley graph} of $G$ with the \textit{connection set} $S$, denoted $\mathrm{Cay}(G, S)$, is a graph whose vertex set is $G$ and two distinct vertices $a$ and $b$ are adjacent if $a^{-1}b \in S$. In particular, if $G = {\mathbb{Z}}_n$ then $\mathrm{Cay}(G, S)$ is called a \textit{circulant graph}. A \textit{cycle} on $n$ vertices, denoted $C_n$, can be regarded as the circulant graph $\mathrm{Cay}({\mathbb{Z}}_n, S)$, where $S = \{1, n-1\}$. The \textit{complement} of the Cayley graph $\mathrm{Cay}(G, S)$ is the Cayley graph $\mathrm{Cay}(G, G \setminus (S \cup \{\textbf{1}\}))$. In particular, the complement of $C_n$, denoted ${\overline{C}}_n$, is  the circulant graph $\mathrm{Cay}({\mathbb{Z}}_n, {\mathbb{Z}}_n \setminus \{0, 1, n-1\})$. 

For a finite group $G$, a \textit{representation} of $G$ is a homomorphism $\rho \colon G \rightarrow \mathrm{GL}(V)$, where $\mathrm{GL}(V)$ is the group of all invertible linear operators on a finite-dimensional complex vector space $V$. The dimension of $V$ is called the \textit{degree} of $\rho$. Two representations $\rho_1 \colon G \rightarrow \mathrm{GL}(V_1)$ and $\rho_2 \colon G \rightarrow \mathrm{GL}(V_2)$ are \textit{equivalent} if there exists an isomorphism $T \colon V_1 \rightarrow V_2$ with the property that $T \rho_1(g) = \rho_2(g)T$ for all $g\in G$. 

The \textit{character} of a representation $\rho$, denoted $\chi_{\rho}$, is the mapping $\chi_\rho \colon G \rightarrow \mathbb{C}$ which is defined as $\chi_\rho(g) = \mathrm{tr}(\rho(g))$ for all $g \in G$, where $\mathrm{tr}(\rho(g))$ denotes the trace of $\rho(g)$. A subspace $U$ of $V$ is said to be $G$-\textit{invariant} if $U$ is invariant under $\rho(g)$ for all $g\in G$. Clearly, $\{\textbf{0}\}$ and $V$ are $G$-\textit{invariant} subspaces of $U$, called the \textit{trivial} subspaces. If $V$ has no non-trivial $G$-invariant subspaces, then $\rho$ is called an \textit{irreducible representation} and $\chi_\rho$ an \textit{irreducible character} of $G$. For an abelian group, each irreducible representation has degree one. Therefore each irreducible representation of an abelian group can be identified with its character. 

In what follows, we write $G = \{a_0, a_1, \ldots, a_{n-1} \}$, where $G$ is an abelian group with the identity element $a_0$. By fundamental theorem of finite abelian groups, we can write
\begin{align*}
G = \mathbb{Z}_{n_1} \oplus \cdots \oplus \mathbb{Z}_{n_k}, 
\end{align*}
where $n_j$ is an integer such that $n_j > 1$ for $1 \leq j \leq k$ and $n = n_1 \cdots n_k$. For $0 \leq r \leq n - 1$, write $a_r = (g_{1r},\ldots,g_{kr})$ and define $$\chi_{r}(a_\ell) = \prod_{j = 1}^{k} {\omega_{n_j}}^{g_{jr} g_{j\ell}}~\mathrm{for~all}~a_\ell  \in G,$$ where $\omega_{n_j} = \exp(\frac{2 \pi \textbf{i}}{n_j})$. Then $\{\chi_r \colon 0 \leq r \leq n - 1\}$ is a complete set of non-equivalent irreducible representations of $G$.


By eigenvalues and eigenvectors of a graph means the eigenvalues and eigenvectors of the adjacency matrix of the graph, respectively. The following lemma is the key ingredient to compute the eigenvalues and the corresponding eigenvectors of an abelian Cayley graph.
\begin{lema}\emph{\cite{Steinberg}}\label{spectra abelian Cayley graph}
Let $S$ be a subset of $G \setminus \{a_0\}$ such that $S = \{-y \colon y \in S\}$. For  $0 \leq r \leq n - 1$, let  $\lambda_r = \sum_{y \in S} \chi_r(y)$ and $\textbf{v}_r = \frac{1}{\sqrt{n}}[\chi_r(x)]_{x \in G}^{T}$. Then 
\begin{enumerate}
\item [(i)] the eigenvalues of $\mathrm{Cay}(G, S)$ are $\lambda_0, \ldots, \lambda_{n-1}$, and 
\item [(ii)] $\textbf{v}_r$ is an eigenvector of $\mathrm{Cay}(G, S)$ corresponding  to the eigenvalue $\lambda_r$ and $\{\textbf{v}_r \colon 0 \leq r \leq n - 1\}$ is an orthonormal basis of eigenvectors.
\end{enumerate} 
\end{lema}
Let $A$ be the adjacency matrix of an abelian Cayley graph $\mathrm{Cay}(G, S)$ such the spectral decomposition of $A$ be given by $A = \sum_{r=0}^{n-1} \lambda_r E_r$, where $E_r = {\textbf{v}}_r {\textbf{v}}^\star_r$. Then the spectral decomposition of the transition matrix is
\begin{align*}
H(t) = \sum_{r=0}^{n-1} \exp(-\textbf{i}t\lambda_r) E_r.
\end{align*}

Let $a_r, a_\ell \in G$ such that $a_r=(g_{1r},\ldots,g_{kr})$ and $a_\ell=(g_{1 \ell},\ldots,g_{k \ell})$. We say that $a_r$ is greater than $a_\ell$, denoted $a_r > a_\ell$, if there exists an $i$ with $1 \leq i \leq k$ such that $g_{jr} = g_{j\ell}$ for $1 \leq j \leq i-1$ and $g_{ir} > g_{i \ell}$. For $\textbf{n}=(n_1,\ldots,n_k)$, let $a_r a_{\ell} = (g_{1r}g_{1 \ell}, \ldots, g_{kr}g_{k \ell})$ and $\frac{a_r}{\textbf{n}} = \left(\frac{g_{1r}}{n_1},\ldots,\frac{g_{kr}}{n_k}\right)$. Further, let
$$\mathrm{wt}(a_r)=\sum_{j=1}^{k}g_{jr},~~\mathrm{wt}(a_r a_\ell)=\sum_{j=1}^{k}g_{jr}g_{j\ell}~~~\mathrm{and}~~~\mathrm{wt}\left(\frac{a_r}{\textbf{n}} \right)=\sum_{j=1}^{k}\frac{g_{jr}}{n_j}.$$

Wang et al.~\cite{J. Wang} provided necessary and sufficient conditions for the existence of FR on abelian Cayley graphs.
\begin{theorem}\label{Fractional Revival Theorem Abelian cayley Graph}\emph{\cite{J. Wang}}
Let $a$ be a non-identity vertex of the abelian Cayley graph $\mathrm{Cay}(G, S)$ such that $a=(z_1,\ldots,z_k)$. Also, let $N$ denote the set
$$\left\{(a_r,a_\ell) \colon a_r,a_\ell \in G, a_r > a_\ell~\mathrm{and}~\mathrm{wt}\left(\frac{2a(a_r-a_\ell)}{\textbf{n}}\right)~\mathrm{is~even}\right\}.$$
Then $\mathrm{Cay}(G, S)$ exhibits FR between $\textbf{0}$ and $a$ if and only if the following three conditions hold.
\begin{enumerate}
\item[(i)] $a$ is an element of order two.
\item[(ii)] $n_j$ is even for $z_j \neq 0$ and $1 \leq j \leq k$.
\item[(iii)] There exists a positive real number $t$ such that $\frac{t}{2\pi}(\lambda_r - \lambda_\ell) \in \mathbb{Z}$ for all $(a_r,a_\ell) \in N$. 
\end{enumerate}
If the conditions $(i),~(ii)$ and $(iii)$ hold and $\lambda_r - \lambda_\ell \in \mathbb{Z}$ for all $(a_r,a_\ell) \in N$, then there is a minimum time $t=\frac{2\pi}{M}$ at which FR occurs between $\textbf{0}$ and $a$, where $M = \gcd (\{\lambda_r - \lambda_\ell \colon (a_r,a_\ell) \in N\})$.
\end{theorem}

Cao and Luo~\cite{Cao} provided another necessary and sufficient conditions for the existence of FR on abelian Cayley graphs.
\begin{theorem}\label{Cao abelian Cayley graph}\emph{\cite{Cao}}
Let $a$ be a non-identity vertex of the abelian Cayley graph $\mathrm{Cay}(G, S)$. Also, let $\alpha$ and $\beta$ be two non-zero complex numbers. Then $\mathrm{Cay}(G, S)$ exhibits ($\alpha,\beta$)-FR between $\textbf{0}$ and $a$ at time $t$ if and only if the following three conditions hold.
\begin{enumerate}
\item[(i)] $a$ is an element of order two.
\item[(ii)] The graph $\mathrm{Cay}(G, S)$ is integral.
\item[(iii)] \[\exp(-\textbf{i}t \lambda_r) = \left\{ \begin{array}{ll} \alpha + \beta & \textrm{ if } r \in X_1\\
 \alpha - \beta & \textrm{ if } r \in X_2, \end{array}\right. \] where $X_1$ and $X_2$ are given by $X_1 = \{r \colon \chi_r(a) = 1\}$ and $X_2 = \{r \colon \chi_r(a) = -1\}$, respectively. 
\end{enumerate}
\end{theorem}

In~\cite{Tan}, Tan et al. proved the following result.
\begin{theorem}\label{Tan abelian Cayley graph}\emph{\cite{Tan}}
Let $\mathrm{Cay}(G, S)$ be an connected integeral abelian Cayley graph. Then $\mathrm{Cay}(G, S)$ is periodic with the minimum period $\frac{2\pi}{R}$, where $R$ is given by $R = \gcd(\lambda_0 - \lambda_r \colon 1 \leq r \leq n - 1)$. 
\end{theorem}

For two real numbers $t$ and $\varepsilon$, we write $|t| < \varepsilon ~(\mathrm{mod} ~2 \pi)$ to mean that $-\varepsilon < t - 2 \pi m < \varepsilon$, for some integer $m$. The following lemma, called \textit{Kronecker approximation theorem}, is useful in  investigating PGFR on abelian Cayley graphs.
\begin{lema}\label{Kronecker}\emph{\cite{Kronecker}}
For an integer $n$ such that $n > 1$, let $\theta_1,\ldots, \theta_{n-1}, \mu_1,\ldots, \mu_{n-1}$ be arbitrary real numbers and $\varepsilon$ be a positive real number. Then the following system of inequalities 
$$|{\theta_r} t - \mu_r| < \varepsilon ~(\mathrm{mod} ~2 \pi) ~\mathrm{for}~r \in \{1, \ldots, n-1\}$$ 
has a simultaneous solution if and only if, for integers $\ell_1,\ldots, \ell_{n-1}$, 
$${\ell_1}{\theta_1} + \cdots + {\ell_{n-1}}{\theta_{n-1}} = 0$$
implies
$${\ell_1}{\mu_1} + \cdots + {\ell_{n-1}}{\mu_{n-1}} \equiv 0 ~(\mathrm{mod}~ 2 \pi).$$
\end{lema}

Let ${\omega}_n = \exp(\frac{2 \pi \textbf{i}}{n})$, where $n$ is a positive integer. The $n$-th \textit{cyclotomic polynomial}, denoted ${\Phi_{n}}(x)$, is defined as 
\begin{align*}
{\Phi_{n}}(x) = \prod_{a \in U(n)}(x - {{\omega}_n}^a),
\end{align*}
where $U(n)=\{r \colon 1 \leq r \leq n-1~\mathrm{and}~\gcd(r,n)=1\}$. Indeed, ${\Phi_{n}}(x)$ $\in {\mathbb{Z}}[x]$ and it is the minimal polynomial for any primitive $n$-th root of unity.  Note that if $n$ is an odd integer with $n > 1$, then ${\Phi_{2n}}(x) = {\Phi_{n}}(-x)$. Also, for an odd prime $p$ and a positive integer $s$, ${\Phi_{p^s}}(x)$ satisfies the relation
\begin{align*}
{\Phi_{p^s}}(x) = \sum_{j=0}^{p-1} x^{j p^{s-1}}.
\end{align*}
For more details on cyclotomic polynomial, see~\cite{Dummit}. Lehmer~\cite{Lehmer} proved the following theorem.
\begin{theorem}\label{cosine rational}\emph{\cite{Lehmer}}
Let $n$ and $k$ be integers such that $n > 2$ and $\gcd(k, n) = 1$. Then $2 \cos(\frac{2 \pi k}{n})$ is an algebraic integer.
\end{theorem}
It is well known that a rational number is an algebraic integer if and only if it is an integer. Therefore, we obtain the following corollary as a consequence of Theorem~\ref{cosine rational}.
\begin{corollary}\label{when does cos is rational}
Let $n$ be a positive integer. Then $\cos(\frac{2 \pi}{n}) \in \mathbb{Q}$ if and only if $n \in \{1, 2, 3, 4, 6\}$.
\end{corollary}

\section{PGFR on abelian Cayley graphs}\label{sec 3}
In this section, we determine a necessary and sufficient condition for the existence of PGFR on abelian Cayley graphs. Let $a$ and $b$ be two distinct vertices of a graph $\Gamma$. If $\textbf{x}$ is an eigenvector of $\Gamma$ such that the entries of $\textbf{x}$ are indexed by the vertices of $\Gamma$, then we use $\widetilde{\textbf{x}}$ to denote the restriction of $\textbf{x}$ to $a$ and $b$. Further, $\textbf{x}_a$ denotes the entry of $\textbf{x}$ at the position corresponding to $a$. The vertices $a$ and $b$ are said to be \textit{strongly cospectral} if for every eigenvector $\textbf{x}$ of $\Gamma$, we have $\textbf{x}_a = \pm \textbf{x}_b$. It follows that $a$ and $b$ are strongly cospectral if and only if for every eigenvector $\textbf{x}$ of $\Gamma$, either $\widetilde{\textbf{x}} = \textbf{0}$ or $\widetilde{\textbf{x}}$ is an eigenvector of the matrix $M$ given by 
$$\begin{pmatrix}
0 & 1  \\
1 & 0
\end{pmatrix}.$$ 

In case of abelian Cayley graph, it is easy to prove that PGFR occurs between two vertices $a$ and $b$ if and only if it occurs between $\textbf{0}$ and $b-a$, where $\textbf{0}$ is the identity element of $G$. Hence, without loss of generality, we consider the existence of PGFR on abelian Cayley graph $\mathrm{Cay}(G, S)$ between the vertices $\textbf{0}$ and $a$ with $a \neq \textbf{0}$. 

Let $\{{\textbf{x}}_0, \ldots, {\textbf{x}}_{n - 1}\}$ be an orthonormal basis of ${\mathbb{C}}^n$ consisting of eigenvectors of $\mathrm{Cay}(G, S)$. Let $\textbf{0}$ and $a$ be two strongly cospectral vertices of $\mathrm{Cay}(G, S)$. We have the assumption that the first two entries of $\textbf{x}_r~(0 \leq r \leq n-1)$ are indexed by the vertices $\textbf{0}$ and $a$. Now we define the following three sets:  
\begin{align*}
X_0 &= \{r \colon 0 \leq r \leq n - 1~\mathrm{and}~\widetilde{{\textbf{x}}}_r = \textbf{0}\},\\
X_1 &= \{r \colon 0 \leq r \leq n - 1~\mathrm{and}~M \widetilde{{\textbf{x}}}_r =  \widetilde{{\textbf{x}}}_r\}~\mathrm{and}\\
X_2 &= \{r \colon 0 \leq r \leq n - 1~\mathrm{and}~M \widetilde{{\textbf{x}}}_r = - \widetilde{{\textbf{x}}}_r\}.
\end{align*}

The following theorem follows from Theorem 2.4 of Chan et al.~\cite{PGFR Adjacency}.   
\begin{theorem}\label{Theorem Chan}\emph{\cite{PGFR Adjacency}}
Let $a$ be a non-identity vertex of the abelian Cayley graph $\mathrm{Cay}(G, S)$. Then $\mathrm{Cay}(G, S)$ exhibits PGFR between $\textbf{0}$ and $a$ if and only if the following two conditions hold.
\begin{enumerate}
\item[(i)] Vertices $\textbf{0}$ and $a$ are strongly cospectral.
\item[(ii)] For arbitrary integers $m_r~(r \in X_1 \cup X_2)$, the following two relations  
$$\sum_{r \in X_1}{m_r} \lambda_r + \sum_{r \in X_2}{m_r} \lambda_r = 0~~\mathrm{and}~~\sum_{r \in X_1 \cup X_2}{m_r} = 0$$
implies that
$$\sum_{r \in X_1}m_{r} \neq \pm 1.$$
\end{enumerate}
\end{theorem}

The following theorem follows from Lemma 4.1 of Arnadottir and Godsil~\cite{Arnadottir}.
\begin{theorem}\label{Strongly Cospectral Vertices}\emph{\cite{Arnadottir}}
Let $a$ be a non-identity vertex of the abelian Cayley graph $\mathrm{Cay}(G, S)$ such that $\textbf{0}$ and $a$ are strongly cospectral. Then the order of $a$ must be two.
\end{theorem}

Recall that the eigenvectors of $\mathrm{Cay}(G, S)$ are given by
$${\textbf{v}}_r = \frac{1}{\sqrt{n}}[\chi_r(\textbf{0}), \chi_r(a), \ldots, \chi_r(a_{n - 1})]^T~~\mathrm{for}~0 \leq r \leq n - 1.$$ 
If $\mathrm{Cay}(G, S)$ exhibits PGFR between $\textbf{0}$ and $a$, then Theorem~\ref{Theorem Chan} and Theorem~\ref{Strongly Cospectral Vertices} implies that $a$ is of order two and $n$ is even. Also, Theorem~\ref{Theorem Chan} gives $X_0 = \emptyset$, $X_1 = \{r \colon \chi_r(a) = 1\}$ and $X_2 = \{r \colon \chi_r(a) = -1\}$.

In the following, we prove a theorem which is a reformulation of Theorem~\ref{Theorem Chan}.
\begin{theorem}\label{PGFR Theorem}
Let $n$ be an even positive integer and $a$ be a vertex of $\mathrm{Cay}(G, S)$ of order two. Then $\mathrm{Cay}(G, S)$ exhibits PGFR between $\textbf{0}$ and $a$ if and only if for arbitrary integers $\ell_1, \ldots, \ell_{n - 1}$, the relation
$$\sum_{r=1}^{n-1}{\ell_r}(\lambda_r - \lambda_0) = 0$$
implies that
$$\sum_{r \in X_2~}\ell_{r} \neq \pm 1.$$
\end{theorem}
\begin{proof}
Let the graph $\mathrm{Cay}(G, S)$ exhibit PGFR between $\textbf{0}$ and $a$, and $\ell_1, \ldots, \ell_{n - 1} \in \mathbb{Z}$ such that 
$$\sum_{r=1}^{n-1}{\ell_r}(\lambda_r - \lambda_0) = 0.$$
This implies that
$$\left (-\sum_{r=1}^{n-1}{\ell_r} \right)\lambda_0 + \sum_{r=1}^{n-1}{\ell_r} \lambda_r = 0.$$
Consider the integers $m_0, m_1, \ldots, m_{n - 1}$ such that
\[m_r = \left\{ \begin{array}{ll} -\sum_{r=1}^{n-1}{\ell_r} & \textrm{ if } r = 0\\
 \ell_r & \textrm{ otherwise. } \end{array}\right. \]
Then $\sum_{r = 0}^{n-1}{m_r} \lambda_r = 0$ and $\sum_{r = 0}^{n-1}{m_r} = 0$. Now Theorem~\ref{Theorem Chan} implies that $\sum_{r \in X_1~}m_{r} \neq \pm 1$. 
This further implies that $\sum_{r \in X_2~}\ell_{r} \neq \pm 1$.

Conversely, assume that for $\ell_1, \ldots, \ell_{n - 1} \in \mathbb{Z}$, the relation
$\sum_{r=1}^{n-1}{\ell_r}(\lambda_r - \lambda_0) = 0$ implies that $\sum_{r \in X_2~}\ell_{r} \neq \pm 1$. Let $m_0, m_1, \ldots, m_{n - 1} \in \mathbb{Z}$ such that 
$$\sum_{r = 0}^{n-1}{m_r} \lambda_r = 0~~\mathrm{and}~~\sum_{r = 0}^{n-1}{m_r} = 0.$$
Then
$$\sum_{r=1}^{n-1}{m_r}(\lambda_r - \lambda_0) = 0.$$
This implies that $\sum_{r \in X_2~}m_{r} \neq \pm 1$. This further implies that $\sum_{r \in X_1~}m_{r} \neq \pm 1$. Therefore, by Theorem~\ref{Theorem Chan} we conclude that the graph $\mathrm{Cay}(G, S)$ exhibits PGFR between $\textbf{0}$ and $a$.
\end{proof}
\begin{rem}\label{PGFR necessary suffficient for circulant}
Let $n$ be an even integer and ${\mathbb{Z}}_n = \{0, \ldots, n-1\} = \{a_{0}, \ldots, a_{n-1}\}$ such that $r$ and $a_{r}$ have the same parity for $1 \leq r \leq n - 1$ and $a_0=0$. Then $\mathrm{Cay}({\mathbb{Z}}_n, S)$ exhibits PGFR between $0$ and $\frac{n}{2}$ if and only if for integers $\ell_1, \ldots, \ell_{n - 1}$, the relation $\sum_{r=1}^{n-1}{\ell_r}(\lambda_r - \lambda_0) = 0$ implies that $\sum_{r~\mathrm{odd}~}\ell_{r} \neq \pm 1$.  
\end{rem}

\section{PGFR on some classes of circulant graphs}\label{sec 4}
Let $n = 2 p^s$, where $p$ is an odd prime and $s \in \mathbb{N}$. In this section, we find a subset $S$ of ${\mathbb{Z}}_n \setminus \{0\}$, consisting of only odd integers, such that $\mathrm{Cay}({\mathbb{Z}}_n, S)$ as well as its complement exhibit PGFR [see Theorem~\ref{PGFR existence example}]. We also find another subset $S$ of ${\mathbb{Z}}_n \setminus \{0\}$, consisting of both odd and even integers, such that $\mathrm{Cay}({\mathbb{Z}}_n, S)$ as well as its complement exhibit PGFR [see Theorem~\ref{PGFR first existence example}]. From these, we obtain two infinite classes of circulant graphs exhibiting PGFR that fails to exhibit FR and PGST. For a suitably chosen $n$, we also prove a necessary condition on the connection set $S$ for the existence of PGFR on $\mathrm{Cay}({\mathbb{Z}}_n, S)$ and its complement [see Theorem~\ref{1st PGFR Nonexistence} and Theorem~\ref{Complement does not exhibit PGFR}]. We also find some classes of circulant graphs that fail to exhibit PGFR see [Corollary~\ref{2pq}, Theorem~\ref{2nd PGFR Nonexistence} and Theorem~\ref{another family not exhibiting PGFR}]. 

Chan et al. proved the following sufficient condition for the existence of PGFR on cycles. 
\begin{lema}\emph{\cite{PGFR Adjacency}}\label{PGFR in cycle}
Let $p$ be an odd prime and $s \in \mathbb{N}$. Then a cycle on $2p^s$ vertices exhibits PGFR.
\end{lema}

We extend Lemma~\ref{PGFR in cycle} to obtain more circulant graphs on $2p^s$ vertices exhibiting PGFR. Recall that we write ${\mathbb{Z}}_n = \{a_{0}, \ldots, a_{n-1}\}$ such that $r$ and $a_{r}$ have the same parity for $1 \leq r \leq n - 1$ and $a_0=0$.
\begin{theorem}\label{PGFR existence example}
Let $n = 2 p^s$, where $p$ is an odd prime and $s \in \mathbb{N}$. Also, let $S$ be a subset of $\mathbb{Z}_n$ such that $S = \{p^{k_0}, \ldots, p^{k_m}, n - p^{k_0}, \ldots, n - p^{k_m}\}$, where $m \in \mathbb{N} \cup \{0\}$ and $0 =k_0 < \cdots < k_m <  s$. If $\gcd(p, m + 1) = 1$, then the circulant graph $\mathrm{Cay}({\mathbb{Z}}_n, S)$ and its complement exhibit PGFR. 
\end{theorem}

\begin{proof}
The eigenvalues of $\mathrm{Cay}({\mathbb{Z}}_n, S)$ are given by
\begin{align*}
\lambda_r = \sum_{i = 0}^{m} ({\omega}_n^{p^{k_i}a_r} + {\omega}_n^{-p^{k_i}a_r}) ~~~~\mathrm{for}~0 \leq r \leq n-1.
\end{align*}
Suppose that the integers $\ell_1,\ldots,\ell_{n - 1}$ satisfy the relation $\sum_{r=1}^{n-1}{\ell_r}(\lambda_r - \lambda_0) = 0$. This implies that
\begin{align*}
\sum_{i = 0}^{m} \sum_{r=1}^{n-1}\ell_r ({\omega}_n^{p^{k_i}a_r} + {\omega}_n^{-p^{k_i}a_r}) - (2m + 2)\sum_{r=1}^{n-1}\ell_r = 0.
\end{align*}
Thus ${\omega}_n$ satisfies the polynomial $f(x)$, where 
\begin{align*}
f(x) = \sum_{i = 0}^{m} \sum_{r=1}^{n-1}\ell_r (x^{p^{k_i}a_r} + x^{n-p^{k_i}a_r}) - (2m + 2)\sum_{r=1}^{n-1}\ell_r.
\end{align*}
Putting $x = -1$, we get  
$$f(-1) = -4(m + 1) \sum_{r~\mathrm{odd}} \ell_r.$$
Since $\Phi_{2 p^s}(x) \in {\mathbb{Z}}[x]$ and it is the minimal polynomial for any primitive $2 p^s$-th root of unity, there exists a polynomial $g(x) \in {\mathbb{Z}}[x]$ such that 
$$f(x) = \Phi_{2 p^s}(x) g(x).$$
Note that $\Phi_{2 p^s}(-1) = \Phi_{p^s}(1) = p$. Therefore we have $$-4(m + 1) \sum_{r~\mathrm{odd}} \ell_r = f(-1) = \Phi_{2 p^s}(-1) g(-1) = pg(-1).$$ 
Since $p$ is an odd prime and $\gcd(p, m + 1) = 1$, the relation $-4(m + 1) \sum_{r~\mathrm{odd}} \ell_r = p g(-1)$ gives that $p$ divides $\sum_{r~\mathrm{odd}} \ell_r$. Thus 
$$\sum_{r ~ \mathrm{odd}~}\ell_{r} \neq \pm 1.$$  
Therefore Remark~\ref{PGFR necessary suffficient for circulant} yields that $\mathrm{Cay}({\mathbb{Z}}_n, S)$ exhibits PGFR. 

The eigenvalues of the complement of $\mathrm{Cay}({\mathbb{Z}}_n, S)$ are given by
\begin{align*}
\lambda_0 = n - 2m - 3~\mathrm{and}~\lambda_r = -1 - \sum_{i = 0}^{m} ({\omega}_n^{p^{k_i}a_r} + {\omega}_n^{-p^{k_i}a_r})~~~~\mathrm{for}~1 \leq r \leq n-1.
\end{align*}
Let $\ell_1,\ldots,\ell_{n - 1}$ be integers such that $\sum_{r=1}^{n-1}{\ell_r}(\lambda_r - \lambda_0) = 0$. This yields that
\begin{align*}
\sum_{i = 0}^{m} \sum_{r=1}^{n-1}\ell_r ({\omega}_n^{p^{k_i}a_r} + {\omega}_n^{-p^{k_i}a_r}) - (2m + 2)\sum_{r=1}^{n-1}\ell_r + n\sum_{r=1}^{n-1}\ell_r = 0.
\end{align*}
Then ${\omega}_n$ satisfies the polynomial given by
\begin{align*}
f(x) = \sum_{i = 0}^{m} \sum_{r=1}^{n-1}\ell_r (x^{p^{k_i}a_r} + x^{n-p^{k_i}a_r}) - (2m + 2)\sum_{r=1}^{n-1}\ell_r + n\sum_{r=1}^{n-1}\ell_r.
\end{align*}
Putting $x = -1$, we get  
$$f(-1) = -4(m + 1) \sum_{r~\mathrm{odd}} \ell_r + n\sum_{r=0}^{n-1}\ell_r.$$
Now the rest of the proof is similar to that of the preceding part, and hence the details are omitted.
\end{proof}

As a consequence of the preceding theorem we obtain the following corollary.
\begin{corollary}\label{Complement cycle exhibiting PGFR}
Let $n = 2 p^s$, where $p$ is an odd prime and $s \in \mathbb{N}$. Then the graphs ${\overline{C}}_n$ exhibits PGFR.
\end{corollary}

The eigenvalues of the graph ${\overline{C}}_n$ are given by $\lambda_0 = n-3$ and $\lambda_r = -1-2\cos(\frac{2 \pi a_r}{n})$ for $1 \leq r \leq n - 1$. Putting $a_r=1$, we obtain $\lambda_r = -1-2\cos(\frac{2 \pi}{n})$. If $n \geq 8$, then Corollary~\ref{when does cos is rational} implies that the graph ${\overline{C}}_n$ is non integral. Therefore, Theorem~\ref{Cao abelian Cayley graph} yields that ${\overline{C}}_n$ does not exhibit FR if $n$ is even and $n \geq 8$. 

It follows from Pal and Bhattacharjya~\cite{Pal 2} [see Corollary 12] that the circulant graphs in Corollary~\ref{Complement cycle exhibiting PGFR} does not exhibit PGST. Thus Corollary~\ref{Complement cycle exhibiting PGFR} provides infinitely many circulant graphs exhibiting PGFR that fail to exhibit PGST and FR.
 
In the following, we obtain more circulant graphs exhibiting PGFR.

\begin{theorem}\label{PGFR first existence example}
Let $n = 2 p^s$, where $p$ is an odd prime and $s \in \mathbb{N}$. Let $y_1,\ldots,y_m \in {\mathbb{Z}}_n$ such that $y_i$ is odd for $1 \leq i \leq \ell$, $y_i$ is even for $\ell+1 \leq i \leq m$ and $\gcd(p, \ell) = 1$. Also, let $S$ is a subset of ${\mathbb{Z}}_n$ such that $S = \{y_i, n - y_i \colon 1 \leq i \leq m \}$. Then the circulant graph $\mathrm{Cay}({\mathbb{Z}}_n, S)$ and its complement exhibit PGFR.
\end{theorem}
\begin{proof}
The eigenvalues of $\mathrm{Cay}({\mathbb{Z}}_n, S)$ are given by
\begin{align*}
\lambda_r = \sum_{i = 1}^{m} ({\omega}_{n}^{y_i a_r} + {\omega}_{n}^{-y_i a_r}) ~~~~\mathrm{for}~0 \leq r \leq n-1.
\end{align*}
Suppose that the integers $\ell_1,\ldots,\ell_{n - 1}$ satisfy the relation $\sum_{r=1}^{n-1}{\ell_r}(\lambda_r - \lambda_0) = 0$. This implies that
\begin{align*}
\sum_{i = 1}^{m} \sum_{r=1}^{n-1}\ell_r ({\omega}_{n}^{y_i a_r} + {\omega}_{n}^{-y_i a_r}) - 2m \sum_{r=1}^{n-1}\ell_r = 0.
\end{align*}
Thus ${\omega}_n$ satisfies the polynomial $f(x)$, where 
\begin{align*}
f(x) = \sum_{i = 1}^{m} \sum_{r=1}^{n-1}\ell_r (x^{y_i a_r} + x^{n-y_i a_r}) - 2m \sum_{r=1}^{n-1}\ell_r.
\end{align*}
Putting $x = -1$, we get  
$$f(-1) = -4\ell \sum_{r~\mathrm{odd}} \ell_r.$$
Then proceeding as in the proof of Theorem~\ref{PGFR existence example}, we find that $\mathrm{Cay}({\mathbb{Z}}_n, S)$ exhibits PGFR. 

The eigenvalues of the complement of $\mathrm{Cay}({\mathbb{Z}}_n, S)$ are given by
\begin{align*}
\lambda_0 = n - 2m - 1~~\mathrm{and}~~\lambda_r = -1 - \sum_{i = 1}^{m} ({\omega}_{n}^{y_i a_r} + {\omega}_{n}^{-y_i a_r})~~~~\mathrm{for}~1 \leq r \leq n-1.
\end{align*}
Let $\ell_1,\ldots,\ell_{n - 1}$ be integers such that $\sum_{r=1}^{n-1}{\ell_r}(\lambda_r - \lambda_0) = 0$. Then we have
\begin{align*}
\sum_{i = 1}^{m} \sum_{r=1}^{n-1}\ell_r ({\omega}_{n}^{y_i a_r} + {\omega}_{n}^{-y_i a_r}) - 2m \sum_{r=1}^{n-1}\ell_r + n\sum_{r=1}^{n-1}\ell_r = 0.
\end{align*}
Therefore ${\omega}_n$ satisfies the polynomial $f(x)$ given by
\begin{align*}
f(x) = \sum_{i = 1}^{m} \sum_{r=1}^{n-1}\ell_r (x^{y_i a_r} + x^{n-y_i a_r}) - 2m \sum_{r=1}^{n-1}\ell_r + n\sum_{r=1}^{n-1}\ell_r.
\end{align*}
Putting $x = -1$, we get  
$$f(-1) = -4\ell  \sum_{r~\mathrm{odd}} \ell_r + n\sum_{r=0}^{n-1}\ell_r.$$
Now the rest of the proof is similar to that of Theorem~\ref{PGFR existence example}.
\end{proof}
The following corollary follows from Theorem~\ref{PGFR first existence example}.
\begin{corollary}\label{PGFR on circulant graphs}
Let $n = 2p^s$, where $p$ is a prime such that $p > 3$ and $s \in \mathbb{N}$. Then the circulant graph $\mathrm{Cay}({\mathbb{Z}}_n, S)$ exhibits PGFR, where $S = \{1, 3, p^s - 3, 2p^s - 1, 2p^s - 3, p^s + 3\}$.
\end{corollary}

Pal proved the following sufficient condition for the non-existence of PGST on circulant graphs.
\begin{theorem}\emph{\cite{Pal}}\label{More circulant graphs}
Let $n = mp$, where $p$ is an odd prime and $m$ is an even positive integer. Also, let $S$ be a subset of ${\mathbb{Z}}_n$ such that  $p \nmid y$ for all $y \in S$. Then the circulant graph $\mathrm{Cay}({\mathbb{Z}}_n, S)$ does not exhibit PGST.
\end{theorem}
The following corollary follows from Theorem~\ref{More circulant graphs}.
\begin{corollary}\label{Not PGST on circulant graphs}
Let $n = 2p^s$, where $p$ is a prime such that $p > 3$ and $s \in \mathbb{N}$. Then the circulant graph $\mathrm{Cay}({\mathbb{Z}}_n, S)$ does not exhibit PGST, where $S = \{1, 3, p^s - 3, 2p^s - 1, 2p^s - 3, p^s + 3\}$.
\end{corollary}
The next lemma tells that the graphs appeared in Corollary~\ref{Not PGST on circulant graphs} fails to exhibit FR.
\begin{lema}\label{Not FR on circulant graphs}
Let $n = 2p^s$, where $p$ is a prime such that $p > 3$ and $s \in \mathbb{N}$. Then the circulant graph $\mathrm{Cay}({\mathbb{Z}}_n, S)$ does not exhibit FR, where $S = \{1, 3, p^s - 3, 2p^s - 1, 2p^s - 3, p^s + 3\}$. 
\end{lema}
\begin{proof}
The eigenvalues of the graph $\mathrm{Cay}({\mathbb{Z}}_n, S)$ are given by
$$\lambda_r = \sum_{y \in S}\cos \left (\frac{ \pi a_r y}{p^s} \right)~\mathrm{for}~0 \leq r \leq n - 1.$$
Putting $a_r = 1$, we obtain $\lambda_r = 2 \cos \left (\frac{\pi}{p^s} \right)$. Since $p > 3$, from Corollary~\ref{when does cos is rational} it is clear that $\mathrm{Cay}({\mathbb{Z}}_n, S)$ is non integral. Now Theorem~\ref{Cao abelian Cayley graph} yields that the graph $\mathrm{Cay}({\mathbb{Z}}_n, S)$ does not exhibit FR.
\end{proof}

From Corollary~\ref{PGFR on circulant graphs}, Corollary~\ref{Not PGST on circulant graphs} and Lemma~\ref{Not FR on circulant graphs}, for a suitably chosen $n$ and a subset $S$ of ${\mathbb{Z}}_n$, we see that $\mathrm{Cay}({\mathbb{Z}}_n, S)$ is an infinite family of circulant graphs exhibiting PGFR that fails to exhibit PGST and FR. Now we present some necessary conditions for some families of circulant graphs exhibiting PGFR. Chan et al.~\cite{PGFR Adjacency} proved the following lemma which gives some cycles not exhibiting PGFR.   
\begin{lema}\emph{\cite{PGFR Adjacency}}\label{PGFR not in cycle} 
A cycle on $n$ vertices does not exhibit PGFR if $n$ is divisible by $2pq$, for some distinct odd primes $p$ and $q$.
\end{lema}
We generalize Lemma~\ref{PGFR not in cycle} in the next theorem. Let $n = mp$, where $m$ is an even positive integer and $p$ is an odd prime. Further, let $p \nmid y$ for all $y \in S$. Pal~\cite{Pal} proved that the eigenvalues of $\mathrm{Cay}({\mathbb{Z}}_n, S)$ satisfy the equation 
\begin{equation}\label{x_14}
(\lambda_2 - \lambda_1) + \sum_{j = 1}^{\frac{p-1}{2}}({\lambda_{mj + 2}} - \lambda_{mj + 1}) + (-1)\sum_{j = 1}^{\frac{p-1}{2}}(\lambda_{mj - 1 } - \lambda_{mj - 2} ) = 0.
\end{equation}
Let $N$ denote a $(n - 2) \times (n - 2)$ unitary matrix and $\Large{\textbf{0}}$ denote the $2 \times (n - 2)$ zero matrix, where $n$ is even. Recall the matrix $M$ from Section~\ref{sec 1} and that   
the first two rows (columns) of $A$ are indexed by $0$ and $\frac{n}{2}$, respectively. 
\begin{theorem}\label{1st PGFR Nonexistence}
Let $n$ be a positive integer such that $2pq \mid n$, where $p$ and $q$ are two distinct odd primes. If the circulant graph $\mathrm{Cay}({\mathbb{Z}}_n, S)$ exhibits PGFR, then there exists $y \in S$ such that either $p \mid y$ or $q \mid y$.
\end{theorem}
\begin{proof}
Let $\mathrm{Cay}({\mathbb{Z}}_n, S)$ exhibit PGFR between the vertices $0$ and $\frac{n}{2}$. Then there exists a sequence of real numbers $\{t_k\}$ such that 
\begin{align*}
\lim_{k\to\infty} H(t_k) =  \begin{pmatrix}
       M & \Large{\textbf{0}}  \\
       {\Large{\textbf{0}}}^T & N
       \end{pmatrix}.
\end{align*}
The first two entries of the eigenvector $\textbf{v}_r$ of $A$ are $\frac{1}{\sqrt{n}}$ and $\frac{1}{\sqrt{n}}{\omega_n}^{a a_r}$, respectively for $0 \leq r \leq n - 1$. Note that $H(t)\textbf{v}_r =  \exp(-\textbf{i}{t}{\lambda_r})\textbf{v}_r$. Since the unit circle is sequentially compact, there exists a subsequence $\{\exp(-\textbf{i}{t^\prime_k}{\lambda_r})\}$ of $\{\exp(-\textbf{i}{t_k}{\lambda_r})\}$ such that $\{\exp(-\textbf{i}{t^\prime_k}{\lambda_r})\}$ converges to a point $\zeta_r$ in the unit circle. Then we have 
\begin{align*}
\lim_{k\to\infty} H(t^\prime_k)\textbf{v}_r = \zeta_r \textbf{v}_r.  
\end{align*}
From the preceding equation, we find that $\zeta_r = \alpha + \beta {\omega_n}^{\frac{n}{2}a_r}$. Therefore, we have  
\begin{equation}\label{x_1}
\lim_{k\to\infty}\exp(-\textbf{i}{t^\prime_k}{\lambda_r}) = \left\{ \begin{array}{cr}
     \alpha + \beta & \mbox{ for even }
     r  \\
     \alpha - \beta &\mbox{ for odd }
     r.
\end{array}\right.
\end{equation}
Assume that $(\alpha + \beta)/(\alpha - \beta) = \exp(-\textbf{i}\eta)$ for some $\eta(\neq 0) \in \mathbb{R}$. Then Equation~(\ref{x_1}) implies that the system of inequalities
\begin{align*}
&|t(\lambda_r - \lambda_{r-1}) - \eta| < \varepsilon~(\mathrm{mod}~ 2\pi)~\mathrm{for~even~r,~and}\\
&|t(\lambda_r - \lambda_{r-1}) + \eta| < \varepsilon~(\mathrm{mod}~2\pi)~\mathrm{for~odd~r}
\end{align*}
has a simultaneous solution for each positive real number $\varepsilon$. 

Note that $n =mp$ for some even positive integer $m$. Suppose that $p \nmid y$ and $q \nmid y$  for all $y \in S$. Then Equation (\ref{x_14}) can be written as 
$${\ell_1}{\theta_1} + \cdots + {\ell_{n-1}}{\theta_{n-1}} = 0,$$
where $\theta_r = \lambda_r - \lambda_{r-1}$ for $1 \leq r \leq n-1$ and 
\[{\ell}_r = \left\{ \begin{array}{rl} 
1 & \textrm{ if } r \in \{mj+2 : 0 \leq j \leq \frac{p-1}{2}\}\\
-1 & \textrm{ if } r \in \{mj-1 : 1 \leq j \leq \frac{p-1}{2}\}\\
  0 & \textrm{ otherwise. } \end{array}\right. \]
Therefore by Kronecker approximation theorem, we find that $p\eta \equiv 0~ (\mathrm{mod}~2\pi)$. Similarly, we also find that $q\eta \equiv 0 ~(\mathrm{mod}~2\pi)$. As $p$ and $q$ are relatively prime, there exist integers $c_1$ and $c_2$ such that $c_1p + c_2q = 1$. Then $\eta = c_1p\eta + c_2q\eta \equiv 0~ (\mathrm{mod}~2\pi)$. This implies that $\beta = 0$, which is a contradiction. Hence, there exists $y \in S$ such that either $p \mid y$ or $q \mid y$.
\end{proof}

Let $n = mp$, where $m$ is an even integer with $m \geq 4$ and $p$ is an odd prime. Further, let $S$ be a connection set in ${\mathbb{Z}}_n$ such that $p \nmid y$ for all $y \in S$. Then Equation~(\ref{x_14}) implies that the eigenvalues of the complement of $\mathrm{Cay}({\mathbb{Z}}_n, S)$ satisfy the equation 
\begin{equation}\label{R_1}
(\lambda_2 - \lambda_1) + \sum_{j = 1}^{\frac{p-1}{2}}({\lambda_{mj + 2}} - \lambda_{mj + 1}) + (-1)\sum_{j = 1}^{\frac{p-1}{2}}(\lambda_{mj - 1 } - \lambda_{mj - 2} ) = 0.
\end{equation}

Using Equation~(\ref{R_1}) and proceeding as in the proof of Theorem~\ref{1st PGFR Nonexistence}, we have the next theorem. 
\begin{theorem}\label{Complement does not exhibit PGFR}
Let $n$ be a positive integer such that $2pq \mid n$, where $p$ and $q$ are two distinct odd primes. If the complement of $\mathrm{Cay}({\mathbb{Z}}_n, S)$ exhibits PGFR, then there exists $y \in S$ such that either $p \mid y$ or $q \mid y$.
\end{theorem}

\begin{corollary}\label{2pq}
Let $n$ be a positive integer such that $2pq \mid n$, where $p$ and $q$ are two distinct odd primes. Then ${\overline{C}}_n$ does not exhibit PGFR.
\end{corollary}

In~\cite{PGFR Adjacency}, Chan et al. proved the following lemma.
\begin{lema}\emph{\cite{PGFR Adjacency}}\label{Chan PGFR Nonexistence 2}
Let $n = 2^{h}p^{s}$, where $p$ is an odd prime, $h, s \in \mathbb{N}$ and $h > 1$. Then the cycle $C_n$ does not exhibit PGFR.
\end{lema}

The following lemma appears in~\cite{PGFR Nonabelian} without proof. We provide a proof of the lemma  borrowing some ideas from~\cite{Bommel}.  

\begin{lema}\emph{\cite{PGFR Nonabelian}}\label{cosine}
Let $p_1$ and $p_2$ be two odd integers such that $1 \leq p_2 < p_1$. Also, let $m = p_1k$, where $k \in \mathbb{N}$. If $d$ is an integer such that $0 \leq d \leq k - 1$, then 
\begin{align*}
\sum_{j = 0}^{p_1 - 1} {(-1)}^{j} \cos \left(\frac{(d + jk)p_2 \pi}{m} \right) = 0.
\end{align*}
\end{lema}
\begin{proof}
We have
\begin{align*}
&\sum_{j = 0}^{p_1 - 1} {(-1)}^{j} \cos \left(\frac{(d + jk)p_2 \pi}{m} \right) \\
=& \cos \left( \frac{dp_2 \pi}{m}\right)\sum_{j = 0}^{p_1 - 1}{(-1)}^{j} \cos \left(\frac{jkp_2 \pi}{m} \right) - \sin \left( \frac{dp_2 \pi}{m}\right)\sum_{j = 0}^{p_1 - 1}{(-1)}^{j} \sin \left(\frac{jkp_2 \pi}{m} \right)\\
=& \cos \left( \frac{dp_2 \pi}{m}\right)\sum_{j = 0}^{p_1 - 1} \cos \left(\frac{j(kp_2+m) \pi}{m} \right) - \sin \left( \frac{dp_2 \pi}{m}\right)\sum_{j = 0}^{p_1 - 1} \sin \left(\frac{j(kp_2+m) \pi}{m} \right).
\end{align*} 
Now  
\begin{align*}
&\sum_{j = 0}^{p_1 - 1} \cos \left(\frac{j(kp_2+m) \pi}{m} \right) + \textbf{i}\sum_{j = 0}^{p_1 - 1} \sin \left(\frac{j(kp_2+m) \pi}{m} \right)\\
=& \sum_{j = 0}^{p_1 - 1} \exp \left(\frac{\textbf{i}j(kp_2+m) \pi}{m} \right)\\
=& \frac{\exp \left(\frac{\textbf{i}p_1(kp_2+m) \pi}{m} \right) - 1}{\exp \left(\frac{\textbf{i}(kp_2+m) \pi}{m} \right) - 1}.
\end{align*} 
Since $p_1$ and $p_2$ are odd, $\frac{p_1(kp_2+m)}{m}$ is even, giving that $\exp \left(\frac{\textbf{i}p_1(kp_2+m) \pi}{m} \right) = 1$. This yields that 
\begin{align*}
\sum_{j = 0}^{p_1 - 1} \cos \left(\frac{j(kp_2+m) \pi}{m} \right) + \textbf{i}\sum_{j = 0}^{p_1 - 1} \sin \left(\frac{j(kp_2+m) \pi}{m} \right) = 0.
\end{align*} 
Thus we have 
$$\sum_{j = 0}^{p_1 - 1} {(-1)}^{j} \cos \left(\frac{(d + jk)p_2 \pi}{m} \right) = 0.$$
This completes the proof.
\end{proof}

Now we generalize Lemma~\ref{Chan PGFR Nonexistence 2} to a larger class of circulant graphs. 
\begin{theorem}\label{2nd PGFR Nonexistence}
Let $n = 2^h p^s$, where $p$ is an odd prime, $h, s \in \mathbb{N}$ and $h > 1$. Also, let $S \subseteq \mathbb{Z}_n$ such that $S = \{p^{k_0}, \ldots, p^{k_m}, n - p^{k_0}, \ldots, n - p^{k_m}\}$, where $m \in \mathbb{N} \cup \{0\}$ and $0 = k_0 < \cdots < k_m < s$. Then $\mathrm{Cay}({\mathbb{Z}}_n, S)$ does not exhibit PGFR. 
\end{theorem}
\begin{proof}
Applying Lemma~\ref{cosine}, we have
\begin{align*}
\sum_{j = 0}^{p^s - 1} {(-1)}^{j} \cos \left(\frac{2(d + 2^{h - 1}j )p^{k_i} \pi}{2^h p^s} \right) = 0,
\end{align*}
for all $i \in \{0, \ldots, m\}$ and $d \in \{0,1,\ldots, 2^{h - 1} - 1\}$. 
This implies that
\begin{equation}\label{x_3}
\sum_{i = 0}^{m}\sum_{j = 0}^{p^s - 1} {(-1)}^{j} \cos \left(\frac{2(d +  2^{h - 1}j)p^{k_i} \pi}{2^h p^s} \right) = 0.
\end{equation}
Note that the eigenvalues $\lambda_r$ of $\mathrm{Cay}({\mathbb{Z}}_n, S)$ are given by
\begin{align*}
\lambda_r = \sum_{i = 0}^{m} 2 \cos \left(\frac{2 r p^{k_i} \pi}{n}\right)~~~~\mathrm{for}~0 \leq r \leq n-1.
\end{align*}
Now we find from Equation (\ref{x_3}) that 
$$ \sum_{j = 0}^{p^s - 1} {(-1)}^{j} (\lambda_{d + 2^{h - 1}j} - \lambda_0) =  -\lambda_0,$$ 
for each $d \in \{0,1, \ldots, 2^{h - 1}-1 \}$. In particular, for $d = 0$ and $d = 1$, we have
\begin{align*}
\sum_{j = 1}^{p^s - 1} {(-1)}^{j} (\lambda_{2^{h - 1}j} - \lambda_0) = - \lambda_0 ~~\mathrm{and}~~
\sum_{j = 0}^{p^s - 1} {(-1)}^{j}( \lambda_{1 + 2^{h - 1}j} - \lambda_0)=  - \lambda_0.
\end{align*}
Consider the integers ${\ell}_r~(1 \leq r \leq n - 1)$ given by
\[{\ell}_r = \left\{ \begin{array}{ll} (-1)^{j} & \textrm{if } r =  2^{h - 1}j,~ \mathrm{where}~1 \leq j \leq p^s - 1\\
 (-1)^{j + 1} & \textrm{if } r = 1 + 2^{h - 1}j,~\mathrm{where}~0 \leq j \leq p^s - 1\\
 0 & \textrm{ otherwise. } \end{array}\right. \]
Then 
$\sum_{r=1}^{n-1}{\ell_r}(\lambda_r - \lambda_0) = 0$ and $\sum_{r~\mathrm{odd}} \ell_r = -1$. Now Remark~\ref{PGFR necessary suffficient for circulant} yields that $\mathrm{Cay}({\mathbb{Z}}_n, S)$ does not exhibit PGFR. 
\end{proof}

In the following we provide a class of complement of cycles not exhibiting PGFR. 
\begin{theorem}\label{another family not exhibiting PGFR}
Let $n = 2^h p^s$, where $p$ is an odd prime, $h, s \in \mathbb{N}$ and $h > 2$. Then the graph ${\overline{C}}_n$ does not exhibit PGFR.  
\end{theorem}
\begin{proof}
The eigenvalues $\lambda_r$ of ${\overline{C}}_n$ are given by
\begin{align*}
\lambda_0=n-3~~~\mathrm{and}~~~\lambda_r = -1- 2 \cos \left(\frac{2 r \pi}{n}\right)~~~~\mathrm{for}~1 \leq r \leq n-1.
\end{align*}
Taking $p_1 = p^s$ and $k = 2^{h-1}$ in Lemma~\ref{cosine}, we obtain 
\begin{align*}
\sum_{j = 0}^{p^s - 1} {(-1)}^{j} \cos \left(\frac{2 (d + 2^{h-1}j)\pi}{n} \right) = 0~~~~\mathrm{for}~0 \leq d \leq 2^{h-1} - 1.
\end{align*}
For $d = 1$ and $d = 2$, we get
$$\sum_{j=0}^{p^s-1}(-1)^j \lambda_{1 +  2^{h-1}j} = -1~~~\mathrm{and}~~~\sum_{j=0}^{p^s-1}(-1)^j \lambda_{2 + 2^{h-1}j} = -1.$$
This further implies that
$$\sum_{j=0}^{p^s-1}(-1)^{j+1} (\lambda_{1 +  2^{h-1}j} - \lambda_0) = 1 + \lambda_0~~~\mathrm{and}~~~\sum_{j=0}^{p^s-1}(-1)^j (\lambda_{2 +  2^{h-1}j} - \lambda_0) = -1 - \lambda_0.$$
Consider the integers given by
\[{\ell}_r = \left\{ \begin{array}{ll} (-1)^{j} & \textrm{if } r = 2 + 2^{h - 1}j,~ \mathrm{where}~0 \leq j \leq p^s - 1\\
 (-1)^{j + 1} & \textrm{if } r = 1 + 2^{h - 1}j,~\mathrm{where}~0 \leq j \leq p^s - 1\\
 0 & \textrm{ otherwise. } \end{array}\right. \]
Then we have $\sum_{r=1}^{n-1}{\ell_r}(\lambda_r - \lambda_0) = 0$ and $\sum_{r~\mathrm{odd}} \ell_r = -1$. Now Remark~\ref{PGFR necessary suffficient for circulant} yields that ${\overline{C}}_n$ does not exhibit PGFR.
\end{proof}
For $n = 4p^s$, where $p$ is an odd prime and $s \in \mathbb{N}$, it is not known whether the graph ${\overline{C}}_n$ exhibits PGFR or not.

\section{PGFR on some classes of non-circulant abelian Cayley graphs}\label{sec 5}
In this section, we explore PGFR on non-circulant abelian Cayley graphs. We present some classes of non-circulant abelian Cayley graphs exhibiting PGFR. We also present some classes of non-circulant abelian Cayley graphs not exhibiting PGFR.

\begin{theorem}\label{more infinite ACG PGFR}
Let $p$ be an odd prime, $s, h \in \mathbb{N}$ such that $s > 1$ and $h > 1$. Let $y_1, \ldots, y_m \in {\mathbb{Z}}_{2p^s}$ and $1 \leq \ell \leq m$ such that $y_i$ is odd for $1 \leq i \leq \ell$, $y_i$ is even for $\ell+1 \leq i \leq m$ and $\gcd(p, \ell) = 1$. Also, let $S$ be a subset of ${\mathbb{Z}}_{2p^s}^h$ such that  $S=\{(y_i,0,\ldots,0),(2p^s-y_i,0,\ldots,0) \colon 1 \leq i \leq m\}$. Then the complement of $\mathrm{Cay}({\mathbb{Z}}_{2p^s}^h, S)$ exhibits PGFR between the vertices $\textbf{0}$ and $a$, where $a=(p^s,0, \ldots, 0)$. 
\end{theorem}
\begin{proof}
It is clear that $X_1 = \{r \colon 0 \leq r \leq 2^h p^{hs}-1, a_r=(g_1, \ldots, g_h), g_i \in {\mathbb{Z}}_{2p^s}~\mathrm{for}~1 \leq i \leq h~\mathrm{and}~g_1~\mathrm{is~even}\}$ and $X_2 = \{r \colon 0 \leq r \leq 2^h p^{hs}-1, a_r=(g_1, \ldots, g_h), g_i \in {\mathbb{Z}}_{2p^s}~\mathrm{for}~1 \leq i \leq h~\mathrm{and}~g_1~\mathrm{is~odd}\}$. 
The eigenvalues of the complement of $\mathrm{Cay}({\mathbb{Z}}_{2p^s}^h, S)$ are given by
$$\lambda_0 = 2^h p^{hs} - 1 - 2m~\mathrm{and}~\lambda_r = -1-\sum_{i=1}^{m}(\omega^{g_1 y_i}_{2p^s} + \omega^{-g_1 y_i}_{2p^s}),~\mathrm{where}~a_r=(g_1, \ldots, g_h)~\mathrm{and}~1 \leq r \leq 2^h p^{hs}-1.$$
Let $\ell_1,\ldots,\ell_{2^hp^{hs} - 1} \in \mathbb{Z}$ such that $\sum_{r=1}^{2^hp^{hs}-1}{\ell_r}(\lambda_r - \lambda_0) = 0$. Then
\begin{align*}
\sum_{i = 1}^{m} \sum_{r=1}^{2^hp^{hs}-1}\ell_r ({\omega}_{2p^s}^{y_i g_1} + {\omega}_{2p^s}^{-y_i g_1}) - 2m \sum_{r=1}^{2^hp^{hs}-1}\ell_r + 2^h p^{hs} \sum_{r=1}^{2^hp^{hs}-1}\ell_r = 0.
\end{align*}
Therefore ${\omega}_{2p^s}$ satisfies the polynomial $f(x)$, where 
\begin{align*}
f(x) = \sum_{i = 1}^{m} \sum_{r=1}^{2^hp^{hs}-1}\ell_r (x^{y_i g_1} + x^{2p^s-y_i g_1}) - 2m \sum_{r=1}^{2^hp^{hs}-1}\ell_r + 2^h p^{hs} \sum_{r=1}^{2^hp^{hs}-1}\ell_r.
\end{align*}
Putting $x = -1$, we get  
$$f(-1) = -4\ell \sum_{r \in X_2} \ell_r + 2^h p^{hs} \sum_{r=1}^{2^hp^{hs}-1}\ell_r.$$
Now following the proof of Theorem~\ref{PGFR existence example}, we obtain that $p \mid \sum_{r \in X_2} \ell_r$ which implies $\sum_{r \in X_2} \ell_r \neq \pm 1$. Therefore,  Theorem~\ref{PGFR Theorem} implies that the complement of $\mathrm{Cay}({\mathbb{Z}}_{2p^s}^h, S)$ exhibits PGFR between the vertices $\textbf{0}$ and $a$.  
\end{proof}

Let $p$ be a prime such that $p > 3$ and $s \in \mathbb{N}$ such that $s > 1$. Consider the complement of the Cayley graph $\mathrm{Cay}({\mathbb{Z}}_{2p^s}^2, S)$, where $S=\{(1,0),(3,0),(p^s-3,0),(2p^s-1,0),(2p^s-3,0),(p^s+3,0)\}$. By Theorem~\ref{more infinite ACG PGFR}, the complement of $\mathrm{Cay}({\mathbb{Z}}_{2p^s}^2, S)$ exhibits PGFR between the vertices $\textbf{0}$ and $a$ such that $a=(p^s,0)$. For $a_r=(1,1)$, the eigenvalue $\lambda_r$ is given by $\lambda_r=-1-2\cos \left(\frac{\pi}{p^s}\right)$. Since $\lambda_r \notin \mathbb{Q}$, therefore Theorem~\ref{Cao abelian Cayley graph} yields that the complement of $\mathrm{Cay}({\mathbb{Z}}_{2p^s}^2, S)$ does not exhibit FR.

In the proof of the next theorem, some ideas are borrowed from  Pal~\cite{Pal}.

\begin{theorem}\label{No PGST on non-circulant}
Let $p$ be an odd prime, $s, h \in \mathbb{N}$ with $s > 1$ and $h > 1$. Let $y_i \in {\mathbb{Z}}_{2p^s}$ such that $p \nmid y_i$ for $1 \leq i \leq m$. Also, let $S$ be a subset of ${\mathbb{Z}}_{2p^s}^h$ such that $S=\{(y_i,0,\ldots,0), (2p^s - y_i,0,\ldots,0) \colon 1 \leq i \leq m \}$. If $a \in {\mathbb{Z}}_{2p^s}^h$ such that $a=(p^s,0, \ldots, 0)$, then the complement of $\mathrm{Cay}({\mathbb{Z}}_{2p^s}^h, S)$ does not exhibit PGST between $\textbf{0}$ and $a$. 
\end{theorem}
\begin{proof}
Assume that the complement of $\mathrm{Cay}({\mathbb{Z}}_{2p^s}^h, S)$ exhibits PGST between $\textbf{0}$ and $a$. Then there exists a sequence of real numbers $\{t_k\}$ and a complex number $\gamma$ of absolute value one such that 
$$\lim_{k\to\infty} H(t_k){\textbf{e}_{\textbf{0}}}=\gamma{\textbf{e}_{a}}.$$
This implies
$$\lim_{k\to\infty} \sum_{r=0}^{2^h p^{hs}-1} \exp(-\textbf{i} t_k \lambda_r) {\chi_r}{(-a)} =  2^h p^{hs}\gamma.$$
Let $a_r=(g_1,\ldots,g_h)$. Since $a=(p^s,0, \ldots, 0)$, we have 
$$\lim_{k\to\infty} \sum_{r=0}^{2^h p^{hs}-1} \exp(-\textbf{i} [t_k \lambda_r + g_1 \pi])=  2^h p^{hs}\gamma.$$
Since the unit circle is sequentially compact, there exists a subsequence $\{\exp(-\textbf{i} [t^\prime_k \lambda_r + g_1 \pi])\}$ of the sequence $\{\exp(-\textbf{i} [t_k \lambda_r + g_1 \pi])\}$ and a complex number $\zeta_r$ of unit modulus for $0 \leq r \leq 2^h p^{hs}-1$ such that  
$$\lim_{k\to\infty} \exp(-\textbf{i}[t^\prime_k \lambda_r + g_1 \pi])=\zeta_r.$$
This implies that
\begin{equation}\label{YY}
\sum_{r=0}^{2^h p^{hs}-1} \zeta_r=2^h p^{hs} \gamma.
\end{equation}
From Equation~(\ref{YY}) we obtain  
$$\zeta_r=\gamma~~~\mathrm{for}~0 \leq r \leq 2^h p^{hs}-1.$$
Therefore we have 
\begin{equation}\label{YYY}
\lim_{k\to\infty} \exp(-\textbf{i}[t^\prime_k \lambda_r+g_1\pi])=\gamma~~~\mathrm{for}~0 \leq r \leq 2^h p^{hs}-1.
\end{equation}
Let $\lambda_g = \lambda_{(g,0,\ldots,0)}$ for $0 \leq g \leq 2p^s-1$. Then Equation~(\ref{YYY}) implies that 
\begin{equation}\label{ADF}
\lim_{k\to\infty} \exp(-\textbf{i}[t^\prime_k \lambda_g+g\pi])=\gamma~~~\mathrm{for}~0 \leq g \leq 2p^s-1.
\end{equation}
This further implies that 
\begin{equation}\label{ADFW}
\lim_{k\to\infty} \exp(-\textbf{i}t^\prime_k [\lambda_{g+1} - \lambda_g])= -1~~~\mathrm{for}~0 \leq g \leq 2p^s-2.
\end{equation}
Let $\lambda^\prime_g$ ($0 \leq g \leq 2p^s-1$) be the eigenvalue of the complement of $\mathrm{Cay}({\mathbb{Z}}_{2p^s}, S^\prime)$, where $S^\prime$ is given by $S^\prime=\{y_1, \ldots, y_m, 2p^s-y_1, \ldots,2p^s - y_m\}$. Then $\lambda^\prime_g = \lambda_g$ for $1 \leq g \leq 2p^s-1$.
Since $p \nmid y_i$ for $1 \leq i \leq m$, Equation~(\ref{R_1}) yields that
\begin{equation}\label{JKLL}
(\lambda_2 - \lambda_1) + \sum_{j = 1}^{\frac{p-1}{2}}({\lambda_{2p^{s-1}j + 2}} - \lambda_{2p^{s-1}j + 1}) + (-1)\sum_{j = 1}^{\frac{p-1}{2}}(\lambda_{2p^{s-1}j - 1 } - \lambda_{2p^{s-1}j - 2} ) = 0.
\end{equation}
Denote the left side of Equation~(\ref{JKLL}) as $L$. Now Equation~(\ref{ADFW}) gives 
$$\lim_{k\to\infty} \exp(-\textbf{i}t^\prime_k L)= -1,$$ which is a contradiction as $L=0$. Thus, the complement of $\mathrm{Cay}({\mathbb{Z}}_{2p^s}^h, S)$ does not exhibit PGST between $\textbf{0}$ and $a$.
\end{proof}
 
The next corollary follows directly from Theorem~\ref{No PGST on non-circulant}.
\begin{corollary} 
Let $p$ be a prime such that $p > 3$ and $s \in \mathbb{N}$ such that $s > 1$. Also, let $S$ be a subset of ${\mathbb{Z}}_{2p^s}^2$ such that $S=\{(1,0),(3,0),(p^s-3,0),(2p^s-1,0),(2p^s-3,0),(p^s+3,0)\}$. Then the complement of $\mathrm{Cay}({\mathbb{Z}}_{2p^s}^2, S)$ does not exhibit PGST between $\textbf{0}$ and $a$, where $a=(p^s, 0)$.
\end{corollary} 

Thus we see that the complement of $\mathrm{Cay}({\mathbb{Z}}_{2p^s}^2, S)$, where ${\mathbb{Z}}_{2p^s}^2$ and $S$ are as mentioned in the preceding corollary, is an infinite family of non-circulant abelian Cayley graphs exhibiting PGFR that fail to exhibit FR and PGST. The next theorem provides a class of graphs exhibiting FR.
\begin{theorem}\label{infinite ACG PGFR}
Let $S \subseteq \mathbb{Z}_2 \oplus {\mathbb{Z}}_p^s$ with $S=\{(1,0,\ldots,0), (0,1,0,\ldots,0),(0,2,0,\ldots,0),\ldots, (0,p-1,0,\ldots,0)\}$, where $p$ is an odd prime and $s > 1$. Then the complement of $\mathrm{Cay}(\mathbb{Z}_2 \oplus {\mathbb{Z}}_p^s, S)$ exhibits FR between the vertices $\textbf{0}$ and $a$ at the minimum time $\frac{2\pi}{p}$, where $a=(1,0, \ldots, 0)$.
\end{theorem}
\begin{proof}
Consider the following four sets
\begin{align*}
X_1^{(1)} &=\{r \colon 0 \leq r \leq 2p^s-1,a_r=(0,g_1, \ldots,g_s),~ g_i \in \mathbb{Z}_p~\mathrm{for}~1 \leq i \leq p~\mathrm{and}~g_1=0\},\\
X_1^{(2)} &=\{r \colon 0 \leq r \leq 2p^s-1,a_r=(0,g_1, \ldots,g_s),~ g_i \in \mathbb{Z}_p~\mathrm{for}~1 \leq i \leq p~\mathrm{and}~g_1 \neq 0\},\\
X_2^{(1)} &=\{r \colon 0 \leq r \leq 2p^s-1,a_r=(1,g_1, \ldots,g_s),~ g_i \in \mathbb{Z}_p~\mathrm{for}~1 \leq i \leq p~\mathrm{and}~g_1=0\}~\mathrm{and}\\
X_2^{(2)} &=\{r \colon 0 \leq r \leq 2p^s-1,a_r=(1,g_1, \ldots,g_s),~ g_i \in \mathbb{Z}_p~\mathrm{for}~1 \leq i \leq p~\mathrm{and}~g_1 \neq 0\}.
\end{align*}
The eigenvalues of the complement of $\mathrm{Cay}(\mathbb{Z}_2 \oplus {\mathbb{Z}}_p^s, S)$ are given by
\[{\lambda}_r=\left\{ \begin{array}{ll} 
2p^s-1-p & \textrm{ if } r=0\\
-1-p & \textrm{ if } r \in X_1^{(1)} \setminus \{0\} \\
-1 & \textrm{ if } r \in X_1^{(2)}\\
1-p & \textrm{ if } r \in X_2^{(1)}\\
1 & \textrm{ if } r \in X_2^{(2)}. \end{array}\right. \]

Recall the set $N$ from Theorem~\ref{Fractional Revival Theorem Abelian cayley Graph}. 
Let $(a_{r_1}, a_{r_2}) \in N$, where $a_{r_1} = (x,x_1,\ldots,x_s)$ and $a_{r_2} = (y,y_1,\ldots,y_s)$. Since $a_{r_1} > a_{r_2}$,  we have $a_{r_1} \neq \textbf{0}$. Also, from the condition that $\mathrm{wt}\left(\frac{2a(a_{r_1}-a_{r_2})}{\textbf{n}}\right)~\mathrm{is~even}$, we find that either $x=y=0$ or $x=y=1$. Thus we have the following eight cases  
 
\noindent \textbf{Case 1.} If $r_1,r_2 \in X_1^{(1)} \setminus \{0\}$ or $r_1,r_2 \in X_1^{(2)}$, then $\lambda_{r_1} - \lambda_{r_2} = 0$.

\noindent \textbf{Case 2.} If $r_1 \in X_1^{(1)} \setminus \{0\}$ and $r_2 \in X_1^{(2)}$, then $\lambda_{r_1} - \lambda_{r_2} = -p$.

\noindent \textbf{Case 3.} If $r_1 \in X_1^{(2)}$ and $r_2 \in X_1^{(1)}\setminus \{0\}$, then $\lambda_{r_1} - \lambda_{r_2} = p$.

\noindent \textbf{Case 4.} If $r_1 \in X_1^{(1)}\setminus \{0\}$ and $r_2 = 0$, then $\lambda_{r_1} - \lambda_{r_2} = -2p^s$.

\noindent \textbf{Case 5.} If $r_1 \in X_1^{(2)}$ and $r_2 = 0$, then $\lambda_{r_1} - \lambda_{r_2} = p-2p^s$.

\noindent \textbf{Case 6.} If $r_1,r_2 \in X_2^{(1)}$ or $r_1,r_2 \in X_2^{(2)}$, then $\lambda_{r_1} - \lambda_{r_2} = 0$.

\noindent \textbf{Case 7.} If $r_1 \in X_2^{(1)}$ and $r_2 \in X_2^{(2)}$, then $\lambda_{r_1} - \lambda_{r_2} = -p$.

\noindent \textbf{Case 8.} If $r_1 \in X_2^{(2)}$ and $r_2 \in X_2^{(1)}$, then $\lambda_{r_1} - \lambda_{r_2} = p$.

The preceding cases implies that $\lambda_{r_1} - \lambda_{r_2} \in \{0,p,-p,-2p^s,p-2p^s\}$ for all $(a_{r_1}, a_{r_2}) \in N$. Observe that $\gcd (\{\lambda_{r_1} - \lambda_{r_2} \colon (a_{r_1},a_{r_2}) \in N\}) = p$. For $t=\frac{2 \pi}{p}$, we obtain $\frac{t}{2 \pi}(\lambda_{r_1} - \lambda_{r_2}) \in \mathbb{Z}$ for all $(a_{r_1}, a_{r_2}) \in N$. Now Theorem~\ref{Fractional Revival Theorem Abelian cayley Graph} implies that the complement of $\mathrm{Cay}(\mathbb{Z}_2 \oplus {\mathbb{Z}}_p^s, S)$ exhibits FR between the vertices $\textbf{0}$ and $a$ at the minimum time $\frac{2\pi}{p}$.  
\end{proof}

It is clear that the graphs appeared in Theorem~\ref{infinite ACG PGFR} are periodic. Applying Theorem~\ref{Tan abelian Cayley graph}, we obtain that the minimum period is $2\pi$. Thus the minimum period of the graphs appeared in the Theorem~\ref{infinite ACG PGFR} is $p$ times the minimum FR time. Now we determine whether the FR obtained in Theorem~\ref{infinite ACG PGFR} is balanced or not. Suppose that the complement of $\mathrm{Cay}(\mathbb{Z}_2 \oplus {\mathbb{Z}}_p^s, S)$ exhibits $(\alpha,\beta)$-FR between $\textbf{0}$ and $a$ at time $\frac{2\pi}{p}$. Then applying Theorem~\ref{Cao abelian Cayley graph}, we obtain $\alpha + \beta = \exp\left(\frac{2\pi \textbf{i}}{p}\right)$ and $\alpha - \beta = \exp\left(-\frac{2\pi \textbf{i}}{p}\right)$. This implies that $\alpha = \cos\left(\frac{2\pi}{p}\right)$ and $\beta = \textbf{i}\sin\left(\frac{2\pi}{p}\right)$. Since $p$ is an odd prime, $|\alpha| \neq |\beta|$. Therefore, the FR is not balanced. 

It follows from Tan et al.~\cite{Tan} [see Theorem 3.5] that the graphs appeared in Theorem~\ref{infinite ACG PGFR} does not exhibit PST. Since the graphs appeared in Theorem~\ref{infinite ACG PGFR} are periodic, Proposition 1.4 of Pal~\cite{State transfer on circulant Pal} implies that they cannot exhibit PGST. Thus we obtain an infinite family of non-circulant abelian Cayley graphs exhibiting FR, and hence exhibiting PGFR, that fails to exhibit PGST. The next theorem provides an infinite family of non-circulant abelian Cayley graphs that fail to exhibit  PGFR.
\begin{theorem}\label{infinite ACG not having PGFR}
Let $m$ be a positive integer such that $m > 3$. Also, let $S$ be a subset of ${{\mathbb{Z}}_2^m}$ such that $S=\{(1,1,1,0,\ldots,0), (0,1,1,0,\ldots,0), (1,1,0,0,\ldots,0)\}$. Further, let $a \in {{\mathbb{Z}}_2^m}$ such that $a=(1, \ldots, 1)$. Then the complement of $\mathrm{Cay}({{\mathbb{Z}}_2^m}, S)$ does not exhibit PGFR between $\textbf{0}$ and $a$. 
\end{theorem}
\begin{proof}
It is clear that $X_1=\{r \colon 0 \leq r \leq 2^m-1~\mathrm{and}~\mathrm{wt}(a_r)~\mathrm{is~even}\}~\mathrm{and}~X_2=\{r \colon 0 \leq r \leq 2^m-1~\mathrm{and}~\mathrm{wt}(a_r)~\mathrm{is~odd}\}$. Consider $r_1 \in X_1$ such that $a_{r_1}=(1,0,1,0,\ldots,0)$ and $r_2 \in X_2$ such that $a_{r_2}=(0,0,1,0,\ldots,0)$. Then $\lambda_{r_1} = \lambda_{r_2}=0$. Define the integers $\ell_1, \ldots, \ell_{2^m-1}$ such that 
\[{\ell}_r=\left\{ \begin{array}{rl} 1 & \textrm{ if } r=r_1\\
 -1 & \textrm{ if } r=r_2~\\
 0 & \textrm{ otherwise. } \end{array}\right. \]
Then $\sum_{r=1}^{2^m-1}{\ell_r}(\lambda_r-\lambda_0)=0$ and $\sum_{r \in X_2} \ell_r=-1$. Now Theorem~\ref{PGFR Theorem} yields that the complement of the graph $\mathrm{Cay}({{\mathbb{Z}}_2^m}, S)$ does not exhibit PGFR between $\textbf{0}$ and $a$.
\end{proof}

\subsection*{Acknowledgements} The first author gratefully acknowledges the support received through the Prime Minister’s Research Fellowship (PMRF), under PMRF-ID: 1903283, funded by the Government of India. The authors also acknowledge the anonymous referees for their valuable comments and suggestions which improve the presentation and quality of the paper.


\end{document}